\documentclass[12pt,a4paper,leqno]{article}
\usepackage{a4}
\usepackage{amsmath}
\usepackage{amssymb}
\newtheorem{theorem}{Theorem}[section]
\newtheorem{lemma}[theorem]{Lemma}
\newtheorem{remark}[theorem]{Remark}
\newtheorem{definition}[theorem]{Definition}

\newcommand{\cc}{\mathbb{C}}
\newcommand{\ccc}{\mathcal{C}}
\newcommand{\dd}{\mathbb{D}}
\newcommand{\mm}{\mathcal{M}}
\newcommand{\nn}{\mathbb{N}}
\newcommand{\oo}{\mathcal{O}}
\newcommand{\qq}{\mathbb{Q}}
\newcommand{\rr}{\mathbb{R}}
\newcommand{\uu}{\mathcal{U}}
\newcommand{\di}{\displaystyle}

\newcommand{\st}{\subset}

\newcommand{\dist}{\operatorname*{dist}}
\newcommand{\rdist}{\operatorname*{r-dist}}

\begin{document}
\title{\bf T-universal Functions With Prescribed Approximation Curves}

\author{D. Mayenberger (Trier)}

\date{}
\maketitle
\begin{abstract}
\noindent Let be $\mathcal{C}$ a family of curves in the unit disc. We show that the set
of all functions $f$ holomorphic on the unit disc, which satisfy the following condition,
is $G_\delta$ and dense in the space of all functions holomorphic on the unit disc.\\
For each compact set $K$ with connected complement, each function $g$ continuous on $K$
and holomorphic on its interior, every point $\zeta_0$ on the unit circle, every curve
$C\in\ccc$ (ending in $\zeta_0$) and any $\varepsilon>0$ there exist numbers $0<a<1$ and
$b\in C$ such that
\[\max_{z\in K}|f(az+b)-g(z)|<\varepsilon \text{ and } |b-\zeta_0|<\varepsilon\]
\end{abstract}

\section{Introduction}
The unit disc will be denoted by $\dd=\{z:|z|<1\}$. The family of compact subsets of
$\cc$ with connected complement we will denote be $\mm$. For any compact set $K\st\cc$ we
will write $A(K)$ for the class of all functions continuous on $K$ and holomorphic in its
interior. The space of all functions holomorphic on $\dd$ endowed with the usual topology
of uniform convergence on compact subsets of $\dd$ will be denoted by $H(\dd)$.\\
In 1976 Luh\cite{luh_1976} proved the existence of a T-universal function $\Phi$ on the
unit disc. T-universality means that translation in the function's $\Phi$ argument forces
the function $\Phi$ to approximate any function $g\in A(K)$ for $K\in\mm$. More
precisely, there exist sequences $0<a_n\to 0$ and $\{b_n\}_n\st \dd$ such that for any
$\zeta\in\partial\dd$, any $K\in\mm$ and any $g\in A(K)$ there exists a strictly
increasing sequence of natural numbers $\{n_k\}_k$
\begin{align*}
&a_{n_k}z+b_{n_k}\to\zeta \quad (k\to\infty) \text{ for all }z\in K\\
&\max_{z\in K}|\Phi(a_{n_k}z+b_{n_k})-g(z)|\to 0 \quad(k\to\infty).
\end{align*}
A similar result has already been proven by Seidel and Walsh\cite{seidel_walsh} in 1941.
For further results the reader is referred to \cite{luh_1979} and
\cite{luh_martirosian_mueller}.\\
A question arising in this context is whether the points $\{b_n\}_n$ above can be chosen
to lie on any curve $C$ belonging to a prescribed family of curves $\ccc$.
Tenthoff\cite{tenthoff} already gave a positive answer on this question by constructing a
function $f$ satisfying the above conditions such that the $b_n$ can be chosen on any
radius $\{z=re^{i\varphi};0\leq r<1\},\varphi\in[0,2\pi)$ of the unit disc. We will proof
this result for general families of curves

\section{T-universal Functions With Prescribed Approximation Curves on the Unit Disc}

\subsection{Continuous Families of Curves}

First we define the notion of a general family of curves. Right after we will restrict
our considerations to those families of curves, which we will call continuous ones.

\begin{definition}
\begin{normalfont}
Let be $I,J\st\rr$ intervals and $z_\alpha:I\longrightarrow\cc$ a continuous function for
each $\alpha\in J$ which satisfies the following conditions
\[
\lim_{t\to\inf(I)}z_\alpha(t)=0 \quad\textnormal{und}
\lim_{t\to\sup(I)}z_\alpha(t)=\infty.
\]
Then we will call the family of functions $(z_\alpha)_{\alpha\in J }$ together with the
intervals $I$ and $J$ a \emph{general family of curves (from zero to infinity)}. Our
short notation will be $\{z_\alpha;I,J\}$.
\end{normalfont}
\end{definition}

\begin{definition}\label{definition_family_ud}
\begin{normalfont}
Let be $I,J\st\rr$ intervals and $z_\alpha:I\longrightarrow\cc$ a continuous function for
each $\alpha\in J$ which satisfies the following conditions
\begin{enumerate}
\item \(\di\lim_{t\to\inf(I)}z_\alpha(t)=z_0\in\dd \quad\textnormal{and}
\lim_{t\to\sup(I)}z_\alpha(t)\in\partial\dd, \) \item for each $\zeta\in\partial\dd$
there exists an $\alpha\in J$ such that $\di\lim_{t\to\sup(I)}z_\alpha(t)=\zeta$.
\end{enumerate}
where $z_0\in\dd$ is fixed and the same for each $\alpha\in J$. Then we will call the
family of functions $(z_\alpha)_{\alpha\in J }$ together with the intervals $I$ and $J$ a
\emph{general family of curves (in the unit disc $\dd$)}. Also here our short notation is
$\{z_\alpha;I,J\}$.
\end{normalfont}
\end{definition}

\begin{definition}
\begin{normalfont}
Let be $I\st\rr$ an interval and $x,y:I\longrightarrow\cc$ two continuous bounded
functions. We set $C_x=x(I)$ und $C_y=y(I)$ and define the \emph{r-distance between $C_x$
and $C_y$} as the number
\[
\rdist(C_x,C_y)=\max_{t\in I}|x(t)-y(t)|.
\]
\end{normalfont}
\end{definition}

\begin{definition}\label{definition_kontinuierlich}
\begin{normalfont}
A general family of curves $\mathcal{C}=\{z_\alpha;I,J\}$ from zero to infinity is to be
called \emph{continuous}, if there will exist a finite or countable subset $\tilde J\st
J$ such that for all $\delta>0,\alpha\in J$ and $j\in\nn$ there exists an
$\tilde\alpha\in\tilde J$ satisfying the following condition:
\[
\rdist\left(z_\alpha(I)\cap\{z:|z|\leq j\}\,,\,z_{\tilde\alpha}(I)\cap\{z:|z|\leq
j\}\right)<\delta.
\]
\end{normalfont}
\end{definition}

Next we give a very simple sufficient criterium for a family of curves to be continuous.

\begin{theorem}
Let be $\mathcal{C}=\{z_\alpha;I,J\}$ a family of curves from zero to infinity such that
the mapping $(\beta,t)\mapsto z_\beta(t)$ is continuous on $J\times I$. Then
$\mathcal{C}$ is a continuous family of curves.
\end{theorem}
\textbf{Proof:}\\
We denote by $D$ the points of $\partial J$ belonging to $J$ and define $\tilde
J=(J\cap\qq)\cup D$. Then $\tilde J$ is countable.\\
Let be given $\delta>0,\alpha\in J,j\in\nn$. Without loss of generality we may assume
that $\delta<1$.\\
For any number $M\in\nn$ we define depending on $I$ the following interval
\begin{equation*}
I_M=\begin{cases} \left(\inf(I),\inf(I)+\frac{1}{M}\right) & \text{, if}
\inf(I)\in(-\infty,\infty),\inf(I)\in I^c \\
(-\infty,-M) & \text{, if} \inf(I)=-\infty \\
\left[\inf(I),\inf(I)+\frac{1}{M}\right) & \text{, if} \inf(I)=\min(I)
\end{cases}
\end{equation*}
By the requirement of the theorem we can choose an $M\in\nn$ such that
\[\sup\{|z_\beta(t)|;t\in I_M,|\beta-\alpha|\leq 1\}<\frac{\delta}{2}.\]
We will fix this $M$. Furthermore we set
\[t_\alpha=\sup\{t\in I:|z_\beta(t)|<2j \text{ for all } \beta \text{ such that }
|\beta-\alpha|\leq 1\}.
\]
Since $z_\alpha(I)$ is a curve from zero to infinity we have
$t_\alpha\in(\inf(I),\sup(I))$. With this number we set
$I_\alpha=[\sup(I_M),t_\alpha]$.\\
Without loss of generality we can assume $\alpha\notin\partial J$ (otherwise we would
have $\alpha\in D\st \tilde J$ and were finished). Thus there exists an $\eta\in(0,1)$
such that $[\alpha-\eta,\alpha+\eta]\st J$.\\
Now $(\beta,t)\mapsto z_\beta(t)$ is continuous and hence uniformly continuous on the
compact set $[\alpha-\eta,\alpha+\eta]\times I_\alpha$. Thus there exists an
$\varepsilon\in(0,1)$, which satisfies the condition
\begin{align*}
&|z_\beta(s)-z_\gamma(t)|<\frac{2\delta}{3}\\
&\text{for all }\beta,\gamma\in[\alpha-\eta,\alpha+\eta],\;s,t\in I_\alpha \text{ with }
|\beta-\gamma|<\varepsilon \text{ and }|s-t|<\varepsilon
\end{align*}
Particularly we have
\[
|z_\alpha(t)-z_\beta(t)|<\frac{2\delta}{3} \text{ for all
}\beta\in(\alpha-\varepsilon,\alpha+\varepsilon),t\in I_\alpha.
\]
Now we choose an $\tilde\alpha\in(\alpha-\varepsilon,\alpha+\varepsilon)$. Then
$\tilde\alpha\in\tilde J$ and we have
\[
\max_{t\in I_\alpha}|z_\alpha(t)-z_{\tilde\alpha}(t)|\leq\frac{2\delta}{3}<\delta
\]
and by the definition of $M$
\begin{equation*}
\max_{t\in I_M}|z_\alpha(t)-z_{\tilde\alpha}(t)|\leq\max_{t\in I_M}|z_\alpha(t)|+
\max_{t\in I_M}|z_{\tilde\alpha}(t)|<\frac{\delta}{2}+\frac{\delta}{2}=\delta.
\end{equation*}
Furthermore by definition of $I_\alpha$ and $I_M$ the following holds
\[
\tilde I=z_\alpha^{-1}\left(\{z:|z|\leq j\}\right)\cup
z_{\tilde\alpha}^{-1}\left(\{z:|z|\leq j\}\right) \st I_M\cup I_\alpha.
\]
By the estimations above we obtain
\[\max_{t\in\tilde I}|z_\alpha(t)-z_{\tilde\alpha}(t)|<\delta\]
and this means
\[
\rdist(z_\alpha(I)\cap\{z:|z|\leq j\}\,,\,z_{\tilde\alpha}(I)\cap\{z:|z|\leq j\})<\delta,
\]
what proves the theorem.\\\\
For the sake of transparency and concreteness we give some examples of continuous
families of curves on the unit disc:
\begin{enumerate}
\item The family of all radii
$\{z_\alpha(t)=te^{i\alpha};t\in[0,1)\}_{\alpha\in[0,2\pi)}$.
\item Logarithmic spirals
$\{z_\alpha(t)=e^{(1+i\alpha)t};t\in(-\infty,\infty)\}_{\alpha\in\rr}$, restricted to the
unit disc.
\item Only one spiral (condition (2) in definition \ref{definition_family_ud} can be
weakened):\\
$z(t)=(1-e^{-t})e^{it},t>0$.
\end{enumerate}

\subsection{Main Result}

First, let be $\ccc=\{z_\alpha;I,J\}$ a fixed continuous family of curves in the unit
disc. Now we are going to define a class of functions we will prove to be nonempty and
moreover a dense $G_\delta$-set in the space $H(\dd)$.\\

The set of all functions $f\in H(\dd)$ such that for every $K\in\mm$, every function
$g\in A(K)$, each $\varepsilon>0$, each $\zeta_0\in\partial\dd$ and any curve
$C\in\mathcal{C}$ (ending in $\zeta_0$) there exist numbers $0<a<1$ and $b\in C$ such
that
\[\max_{z\in K}|f(az+b)-g(z)|<\varepsilon,\quad |b-\zeta_0|<\varepsilon\]
will be called the \emph{class of T-universal functions with respect to the family of
curves
$\ccc$} in $\dd$. It is denoted by $\uu_\ccc(\dd)$.\\
Note that the definition above already implies $aK+b\st\dd$, otherwise the function $f$
would not be defined on $aK+b$.\\

\begin{theorem}\label{main_theorem}
The set $U_\ccc(\dd)$ of T-universal functions on the unit disc with prescribed
approximation curves is $G_\delta$ and dense in $H(\dd)$.
\end{theorem}

\subsubsection{Proof of Theorem \ref{main_theorem}}

We will fix some sequences and sets for abbreviation purposes.

\begin{enumerate}
\item For each $m\in\nn$ we denote $L_m=\{z:|z|\leq m\}$.
\item Let be $\{p_j\}_j$ an enumeration of all polynomials with coefficients in $\qq+ i\qq$.
\item The sequence $\{\zeta_p\}_p$ is chosen to be dense on the unit circle $\partial\dd$.
\item For each $p\in\nn$ we choose with respect to those curves ending in $\zeta_p$ a sequence of
curves $\{C_{pl}\}_l$ according to definition \ref{definition_kontinuierlich}. I.e. for
each curve $C\in\ccc$ ending in $\zeta_p$ we find an index $l\in\nn$ such that $C_{pl}$
lies arbitrarily near to $C$ in terms of the r-distance.
\item For $p,l\in\nn$ we choose a sequence of points $\{b_{nlp}\}_n$ being dense on $C_{pl}$.
\item The sequence $\{a_k\}_k$ is a sequence of positive numbers dense in $(0,1)$.
\end{enumerate}

With this notions we will prove three technical lemmas.  For an intermediate step we need
an auxiliary class. For this we fix an $h\in\nn$. The set of all functions $f\in H(\dd)$
such that for every $K\in\mm(\cc)$, every function $g\in A(K)$, each $\varepsilon>0$,
each $\zeta_0\in\partial\dd$ and any curve $C\in\mathcal{C}$ (ending in $\zeta_0$) there
exist numbers $0<a<1$ and $b\in U_{\frac{1}{h}}(C)=\{z\in\cc:\dist(z,C)<\frac{1}{h}\}$
such that
\[\max_{z\in K}|f(az+b)-g(z)|<\varepsilon,\quad |b-\zeta_0|<\varepsilon\]
will be denoted by $\uu^{(h)}_\ccc(\dd)$.\\

For $m,j,p,s,t,l,k,n\in\nn$ we set
\begin{align*}
\oo_{\mathcal{C}}(m,j,p,s,t,l,k,n)=\Bigl\{&g\in H(\dd): \max_{z\in L_m}
|g(a_kz+b_{nlp})-p_j(z)|<\frac{1}{s};\\
&b_{nlp}\in C_{lp},\;|b_{nlp}-\zeta_p|<\frac{1}{t}\Bigr\}
\end{align*}
Note that this set depends on $\ccc$, although it does not appear itself in the above
definition. But the $C_{pl}$ are chosen in $\ccc$.\\

Our first lemma states that the class $\uu_\ccc(\dd)$ has a representation with
intersections and unions of the sets $\oo_{\mathcal{C}}(m,j,p,s,t,l,k,n)$.

\begin{lemma}\label{darstellung_d}
The following equations hold
\[
\uu_\ccc(\dd)=\bigcap_{h=1}^\infty \uu^{(h)}_\ccc(\dd)=\bigcap_{m,j,p,s,t,l=1}^\infty
\bigcup_{k,n=1}^\infty \oo_{\mathcal{C}}(m,j,p,s,t,l,k,n).
\]
\end{lemma}
\textbf{Proof:}\\
The first equation is obvious due to the definition of the considered classes.\\
Let be $f$ an element of the right hand side an let be given $K\in\mm, g\in A(K),$
$\varepsilon>0,\zeta_0\in\partial\dd$ and a curve $C\in\ccc$ ending in $\zeta_0$. We
fix an $h\in\nn$.\\
Then there is an $m\in\nn$ such that $K\st L_m$. Furthermore we find $s,t\in\nn$ such
that $\frac{1}{s}<\frac{\varepsilon}{2},\frac{1}{t}<\frac{\varepsilon}{2}$. Then by
Mergelyan's theorem we choose a $j\in\nn$ satisfying
\[
\max_K|p_j(z)-g(z)|<\frac{\varepsilon}{2}.
\]
Since the sequence $\{\zeta_p\}_s$ is dense in $\partial\dd$ there is a $p\in\nn$ with
$|\zeta_p-\zeta_0|<\frac{\varepsilon}{2}$.\\
The family of curves $\ccc$ is continuous, so we can find an $l\in\nn$ to this $p$
satisfying
\[
\rdist(C,C_{pl})<\frac{1}{h}.
\]
Due to the definition of $\oo_{\mathcal{C}}(m,j,p,s,t,l,k,n)$ and the representation of
the right hand side there exist numbers $n,k\in\nn$ with the following properties
\[
\max_{L_m}|f(a_kz+b_{nlp})-p_j(z)|<\frac{1}{s}\quad\textnormal{and}\quad b_{nlp}\in
C_{pl},\,|b_{nlp}-\zeta_p|<\frac{1}{t}.
\]
Hence we obtain
\begin{align*}
&\max_K|f(a_kz+b_{nlp})-g(z)|\leq\\
&\max_{L_m}|f(a_kz+b_{nlp})-p_j(z)|+\max_K|p_j(z)-g(z)|<\frac{1}{s}+\frac{\varepsilon}{2}<\varepsilon
\end{align*}
Since $\rdist(C,C_{pl})<\frac{1}{h}$ and $b_{nlp}\in C_{pl}$ it is also true that
$b_{nlp}\in U_\frac{1}{h}(C)$ and hence $f\in\uu_\ccc^{(h)}(\dd)$. Since $h\in\nn$
was arbitrary we conclude $f\in\uu_\ccc(\dd)$.\\
Now let be $f$ a function lying in $\uu_\ccc(\dd)$ and let be given\\
$m,j,p,s,t,l\in\nn$. By definition of $\uu_\ccc(\dd)$ there exist $0<a<1$ and $b\in
C_{pl}$ satisfying
\[
\max_{L_m}|f(az+b)-p_j(z)|<\frac{1}{2s}\quad\textnormal{und}\quad
|b-\zeta_p|<\frac{1}{2t}.
\]
If we set $d=\dist(aL_m+b,\partial\dd)>0$ and
\[
\tilde L_m=\left\{z\in\cc:\dist(az+b,\partial\dd)\geq\frac{d}{2},az+b\in\dd\right\},
\]
then $a\tilde L_m+b$ will be a compact subset of $\dd$ with $a\tilde L_m+b\supset
aL_m+b$. Since $f$ is uniformly continuous on this compact set, there exists a $\delta>0$
such that $|f(z_1)-f(z_2)|<\frac{1}{2s}$
for all $z_1,z_2\in a\tilde L_m+b$, $|z_1-z_2|<\delta$.\\
Then we find numbers $k,n\in\nn$ with $|a_k-a|<\frac{\delta}{2m}$ and
$|b_{nlp}-b|<\min\left\{\frac{\delta}{2},\frac{1}{2t}\right\}$. Thus for all $z\in L_m$
we have $|a_kz+b_{nlp}-(az+b)|<\delta$ and hence we obtain
\begin{align*}
&\max_{L_m}|f(a_kz+b_{nlp})-p_j(z)|\leq\\
&\max_{L_m}|f(a_kz+b_{nlp})-f(az+b)|+\max_{L_m}|f(az+b)-p_j(z)|<\frac{1}{2s}+\frac{1}{2s}=
\frac{1}{s}
\end{align*}
Finally we have
$|b_{nlp}-\zeta_p|\leq|b_{nlp}-b|+|b-\zeta_p|<\frac{1}{2t}+\frac{1}{2t}=\frac{1}{t}$ and
hence $f$ lies in the right hand side of the stated equation. This
proves the lemma.\\\\

The next lemma states, taken together with the preceding one, that $\uu_\ccc(\dd)$ is
indeed a $G_\delta$ set in $H(\dd)$

\begin{lemma}\label{menge_offen_d}
For all $m,j,p,s,t,l,k,n\in\nn$ the set $\oo_\ccc(m,j,p,s,t,l,k,n)$ is open in $H(\dd)$.
\end{lemma}
\textbf{Proof:}\\
Fix $m,j,p,s,t,l,k,n\in\nn$ and $f\in\oo_\ccc(m,j,p,s,t,l,k,n)$. We set
\[\delta=\frac{1}{s}-\max_{L_m}|f(a_kz+b_{nlp})-p_j(z)|>0\]
and define
\[U_\delta(f)=\left\{g\in H(\dd):\max_{L_m}|g(a_kz+b_{nlp})-f(a_kz+b_{nlp})|<\delta\right\}.\]
Then we obtain for all $g\in U_\delta(f)$:
\begin{align*}
&\max_{L_m}|g(a_kz+b_{nlp})-p_j(z)|\leq\\
&\max_{L_m}|g(a_kz+b_{nlp})-f(a_kz+b_{nlp})|+\max_{L_m}|f(a_kz+b_{nlp})-p_j(z)|<\\
&\frac{1}{s}-\max_{L_m}|f(a_kz+b_{nlp})-p_j(z)|+\max_{L_m}|f(a_kz+b_{nlp})-p_j(z)|=
\frac{1}{s}
\end{align*}
Thus the open $\delta$-neighborhood $U_\delta(f)$ of $f$ is contained in
$\oo_\ccc(m,j,p,s,t,l,k,n)$ and the statement follows.\\

\begin{lemma}\label{menge_dicht_d}
For all $m,j,p,s,t,l\in\nn$ the set
\[\bigcup_{k,n=1}^\infty \oo_\ccc(m,j,p,s,t,l,k,n)\]
is dense in the space $H(\dd)$.
\end{lemma}
\textbf{Proof:}\\
Fix numbers $m,j,p,s,t,l\in\nn$ and let be given $f\in H(\dd)$, a compact set $K\st\dd$
and an
$\varepsilon>0$.\\
First we find a compact set $K\st B\st\dd$ with connected complement in $\dd$. Since $B$
is a compact subset of \ $\dd$ we can choose a $\delta>0$ such that\\
$\{z:|z-\zeta_p|<\delta\}\cap B=\varnothing$. Since $C_{pl}$ is a curve ending on the
unit circle and $\{a_k\}_k$ is dense in the interval $(0,1)$ we can choose numbers
$k,n\in\nn$ such that $|b_{nlp}-\zeta_p|<\min\left\{\frac{1}{t},\frac{\delta}{2}\right\}$
and $0<a_k<\frac{\delta}{2m}$. By definition the point $b_{nlp}$ lies in $C_{pl}$
anyway.\\
Thus we have for all $z\in L_m$:
\[
|a_kz+b_{nlp}-\zeta_p|\leq a_k|z|+|b_{nlp}-\zeta_p|<
\frac{\delta}{2m}m+\frac{\delta}{2}=\delta
\]
and hence
\[a_kL_m+b_{nlp}\st \{z:|z-\zeta_p|<\delta\} \st B^c.\]
Due to Runge's theorem on polynomial approximation there exists a polynomial $p$ with
\begin{align*}
&\max_{B}|p(z)-f(z)|<\varepsilon\\
&\max_{L_m}|p(a_kz+b_{nlp})-p_j(z)|<\frac{1}{s}
\end{align*}
Thus $p$ lies "near to" $f$ and we have $p\in\oo_\ccc(m,j,p,s,t,l,k,n)$, what completes
the proof.\\

Now theorem \ref{main_theorem} is a consequence of the preceding three lemmas. Indeed,
lemma \ref{darstellung_d} and \ref{menge_offen_d} state that $\uu_\ccc(\dd)$ is a
$G_\delta$ set in $H(\dd)$. Together with lemma \ref{menge_dicht_d} we obtain that
$\uu_\ccc(\dd)$ has an representation as a countable intersection of dense sets.
Recalling that $H(\dd)$ is a complete metric space and applying Baire's category theorem
we obtain that $\uu_\ccc(\dd)$ is dense in $H(\dd)$. This proves theorem
\ref{main_theorem}.

\subsubsection{An additional property of T-universal functions with prescribed approximation curves}

\begin{remark}
Every function $f\in H(\dd)$ can be expressed as the sum of two T-universal functions
with prescribed approximation curves.
\end{remark}
The following short \textbf{proof} is due to J.-P. Kahane \cite{kahane}.\\
Let be given a function $f\in H(\dd)$. The mapping
\[T_f(g):H(\dd)\to H(\dd),\;T_f(g)=g+f\quad (g\in H(\dd))\]
is a homeomorphism.\\
Since the set $\uu_\ccc(\dd)$  is $G_\delta$ and dense in $H(\dd)$, the same holds for
\[T_f(\uu_\ccc(\dd))=\uu_\ccc(\dd)+f.\]
Due to Baire's category theorem we have
\[\uu_\ccc(\dd)\cap(\uu_\ccc(\dd)+f)\neq\varnothing.\]
Thus there exist $g,h\in\uu_\ccc(\dd)$ with $f=g-h$. Since $-h\in\uu_\ccc(\dd)$ the
result follows.

\end{document}